\documentclass[11pt]{article}
\usepackage{amssymb}
\usepackage{amsmath}

\newtheorem{theorem}{Theorem}
\newtheorem{lemma}{Lemma}

\begin{document}

\title{Explicit Solutions of the One-dimensional Vlasov-Poisson System with Infinite Mass}
\author{Stephen Pankavich \\
Department of Mathematics \\
Indiana University \\
Bloomington, IN 47401 \\
sdp@indiana.edu}
\date{\today}
\maketitle

\begin{center}
\end{center}

\begin{abstract}
A collisionless plasma is modelled by the Vlasov-Poisson system in
one-dimension. A fixed background of positive charge, dependent
only upon velocity, is assumed and the situation in which the
mobile negative ions balance the positive charge as $\vert x \vert
\rightarrow \infty$ is considered.  Thus, the total positive
charge and the total negative charge are infinite.  In this paper,
the charge density of the system is shown to be compactly
supported.  More importantly, both the electric field and the
number density are determined explicitly for large values of
$\vert x \vert$.
\end{abstract}

\section*{Introduction}

Consider the one dimensional Vlasov-Poisson system with a given
positive background function and initial data. We take as given the
functions $F : \mathbb{R} \rightarrow [0,\infty)$ and $f_0 :
\mathbb{R}^2 \rightarrow [0,\infty)$ and seek functions $f : [0,T]
\times \mathbb{R} \times \mathbb{R} \rightarrow \mathbb{R}$ and $E :
[0,T] \times \mathbb{R} \rightarrow \mathbb{R}$ such that

\begin{equation}
\label{VP} \left. \begin{array}{ccc}
& & \partial_t f + v \ \partial_x f - E \ \partial_v f = 0,\\
\\
& & \rho(t,x) = \int \left ( F(v) - f(t,x,v) \right ) \ dv, \\
\\
& & E(t,x) = \frac{1}{2} \left ( \int_{-\infty}^x \rho(t,y) \ dy -
\int_x^\infty \rho(t,y) \ dy \right ),\\
\\
& & f(0,x,v) = f_0(x,v). \end{array} \right \}
\end{equation}

\noindent Here $t \in [0,T]$ denotes time, $x \in \mathbb{R}$
denotes one-dimensional space, and $v \in \mathbb{R}$ denotes
one-dimensional momentum. In (\ref{VP}), $F$ represents a number
density in phase space of positive ions which form a fixed
background, and $f$ describes the number density of mobile negative
ions.  Notice that if $f_0 = F$, then $f = F$ is a steady solution.
Thus, we consider solutions for which $f \rightarrow F$ as $\vert x
\vert \rightarrow \infty$.
 The existence of a unique, local-in-time solution to (\ref{VP}) was shown
in \cite{Pankavich3}.  Thus, we will include assumptions which
satisfy the requirements of \cite{Pankavich3} so that the existence
of a solution $f$ on $[0,T]$ is valid for some $T > 0$, and the
results which follow may be applied to the local-in-time solution.\\

The Vlasov-Poisson system has been studied extensively in the case
where $F(v) = 0$ and solutions tend to zero as $\vert x \vert
\rightarrow \infty$, both for the one-dimensional problem and the
more difficult, three-dimensional problem.  Most of the literature
involving the one-dimensional Vlasov-Poisson system focus on time
asymptotics, such as \cite{BKR} and \cite{BFFM}.  Much more work has
been done concerning the three-dimensional problem.  Smooth
solutions were shown to exist globally-in-time in \cite{Pfaf} and
independently in \cite{LP}.  The results of \cite{Pfaf} were later
revised in \cite{Sch}.  Important results preliminary to the
discovery of a global-in-time solution include \cite{Batt} and
\cite{Horst}. A complete discussion of the literature concerning the
Vlasov-Poisson system may be found in both \cite{Glassey} and \cite{Rein}.\\

Only over the past decade has some work begun studying solutions of
the Vlasov-Poisson system with infinite mass and energy. Under
differing assumptions, distributional solutions with infinite mass
or infinite kinetic energy have been constructed in \cite{Jabin} and
\cite{Perthame}. More recently, the three-dimensional analogue of
(\ref{VP}), which yields solutions with both infinite mass and
energy, has been studied. Local existence of smooth solutions and a
continuation criteria for this problem were shown in \cite{VPSSA}. A
priori bounds on the current density were achieved in \cite{SSAVP}.
Finally, global existence in the case of a radial electric field was
shown in \cite{Pankavich}, and global existence without the
assumption of radial symmetry, in \cite{Pankavich2}.

\section*{Section 1}

We will make the following assumptions throughout, for any $x,v
\in \mathbb{R}$ :

\begin{enumerate}
\renewcommand{\labelenumi}{(\Roman{enumi})}
\item There is $R > 0$ such that for $\vert x \vert
> R$,
\begin{equation}
\label{IC} f_0(x,v) = F(v)
\end{equation}
where $F  \in \mathcal{C}_c^1(\mathbb{R})$ is even and nonnegative,
and $f_0 \in \mathcal{C}^1_c(\mathbb{R}^2)$ is nonnegative.

\item There is $C^{(1)}
> 0 $ such that
\begin{equation}
\label{star} \int_0^T \Vert E(\tau) \Vert_\infty \ d\tau \leq
C^{(1)}.
\end{equation}
\end{enumerate}

Notice that assumption $(I)$ satisfies the assumptions made in
\cite{Pankavich3} and thus the existence of a unique local-in-time
solution is guaranteed.  Define the characteristics $X(s,t,x,v)$
and $V(s,t,x,v)$ by
\begin{equation}
\label{char} \left. \begin{array}{ccc} & &
\frac{\partial}{\partial s} X(s,t,x,v) = V(s,t,x,v) \\
& & \frac{\partial}{\partial s} V(s,t,x,v) = - E(s,X(s,t,x,v)) \\
& & X(t,t,x,v) = x \\
& & V(t,t,x,v) = v \end{array} \right \}
\end{equation}
Then,
\begin{eqnarray*}
\frac{\partial}{\partial s} f(s,X(s,t,x,v),V(s,t,x,v)) & = &
\partial_t f(s,X(s,t,x,v),V(s,t,x,v)) \\
& & + V(s,t,x,v) \partial_x
f(s,X(s,t,x,v),V(s,t,x,v)) \\
& & - E(s,X(s,t,x,v)) \partial_v f(s,X(s,t,x,v),V(s,t,x,v)) = 0
\end{eqnarray*}
so that $f$ is constant along characteristics, and

\begin{equation}
\label{fchar}
f(t,x,v) = f(0, X(0,t,x,v), V(0,t,x,v)) =
f_0(X(0,t,x,v),V(0,t,x,v)).
\end{equation}
Thus, we find that $f$ must be nonnegative and $ \Vert f
\Vert_\infty = \Vert f_0 \Vert_\infty < \infty$.\\

Now, put $$ g(t,x,v) = F(v) - f(t,x,v)$$ and define for every $t
\in [0,T]$,
\begin{equation}
\label{Qg} Q_{g}(t) : = \sup \{ \vert v \vert : \exists \ x \in
\mathbb{R}, \tau \in [0,t] \ \mathrm{such \ that} \ g(\tau,x,v)
\neq 0 \},
\end{equation}
and
\begin{equation}
\label{R}
 R(t) := R + tQ_g(t) + \int_0^t \int_\tau^t \Vert E(s)
\Vert_\infty \ ds \ d\tau.
\end{equation}
Then, $V(s,t,x,v)$ is linear in $v$  for large $\vert x \vert$,
and $\rho(t,x)$ has compact support for every $t \in [0,T]$.

\begin{theorem}
Let $T > 0$ and $f$, a $C^1$ solution of (\ref{VP}) on $[0,T]
\times \mathbb{R}^2$, be given and assume (\ref{IC}) and
(\ref{star}) hold. Then, for $t \in [0,T]$ and $\vert x \vert >
R(t)$ we have

\begin{enumerate}

\item $\rho(t,x) = 0.$

\item For $s \in [0,t]$ and $\vert v \vert \leq Q_g(t)$,
$$ \frac{dV}{dv}(s,t,x,v) = 1.$$
\end{enumerate}

\end{theorem}
More importantly, we may explicitly determine the unique solution
to (\ref{VP}) for large $\vert x \vert$.  Define $$ sign(x) :=
\left \{ \begin{array}{rl} 1, & x \geq 0 \\ -1, & x < 0
\end{array} \right.$$

\begin{theorem}
Let $T > 0$ and $f$, a $C^1$ solution of (\ref{VP}) on $[0,T]
\times \mathbb{R}^2$, be given and assume (\ref{IC}) and
(\ref{star}) hold. Then, for $t \in [0,T]$ and $\vert x \vert >
R(t)$ we have
$$ E(t,x) = E^0 sign(x) \cos (\omega t) $$
and
$$ f(t,x,v) = F\left ( v + \frac{E^0 sign(x)}{\omega} \sin (\omega t) \right ),$$
 where
$$ E^0 = \frac{1}{2} \int_{-R}^R \int (F(v) - f_0(y,v)) \ dv \ dy $$
and
$$ \omega = \left ( \int F(v) \ dv \right )^\frac{1}{2}.$$
\end{theorem}

We derive the form of $f(t,x,v)$ and $E(t,x)$ from the field bound
and the initial current density. Then, using Theorem $1$, we show
that the current density must be a multiple of the field for large
$\vert x \vert$. In order to arrive at these results, we must
first show that $Q_g$ is bounded.

\begin{lemma}
For any $t \in [0,T]$, $$ Q_g(t) \leq C.$$
\end{lemma}
\vspace{0.25in}
\noindent To control spatial characteristics, we
will use the following lemma :

\begin{lemma}
For $t \in [0,T]$, $s \in [0,t]$, $x \in \mathbb{R}$, and $\vert v
\vert \leq Q_g(t)$, we have
$$\vert x \vert \geq R(t) \Rightarrow \vert X(s,t,x,v) \vert \geq
R(s).$$
\end{lemma}

\vspace{0.25in} We will delay the proofs of the lemmas until
Section $4$. Section $2$ will be dedicated to proving Theorem $1$
using the above lemmas. Then, in Section $3$, we will prove
Theorem $2$, utilizing the first theorem.

\section*{Section 2}
In order to prove Theorem $1$, we must first bound the charge
density and $v$-derivatives of characteristics. Notice from
(\ref{IC}) and Lemma $2$,
$$ \vert x \vert > R(t) \Rightarrow f_0(X(s,t,x,v),V(s,t,x,v)) =
F(V(s,t,x,v))$$ for every $s \in [0,t]$, $t \in [0,T]$, and $\vert
v \vert \leq Q_g(t)$.  In particular,
\begin{equation}
\label{ftoF}
\vert x \vert > R(t) \Rightarrow f(t,x,v) =
F(V(0,t,x,v)).
\end{equation}
Let $\vert x \vert > R(t)$ and using (\ref{star}), (\ref{ftoF}),
and Lemma $1$, we write
\begin{eqnarray*}
\vert \rho(t,x) \vert & = & \left \vert \int (F(v) -
f(t,x,v)) \ dv \right \vert \\
& = & \left \vert \int_{\vert v \vert \leq Q_g(t)} (F(v) - F(V(0,t,x,v)) \ dv \right \vert \\
& \leq & \int_{\vert v \vert \leq Q_g(t)} \Vert F^\prime
\Vert_\infty \vert
V(0,t,x,v) - v \vert \ dv \\
& \leq & Q_g(t) \ \Vert F^\prime \Vert_\infty \left (
\int_0^t \Vert E(\tau) \Vert_\infty \ d\tau \right ) \\
& \leq & C
\end{eqnarray*}
So, for $\vert x \vert > R(t)$, \begin{equation} \label{rhobound}
\rho(t,x) \leq C.
\end{equation}

To bound derivatives of characteristics, we use (\ref{VP}) and
(\ref{char}) to find
$$ \frac{\partial \dot{X}}{\partial v}(s,t,x,v) =
\frac{\partial V}{\partial v}(s,t,x,v) $$ and
\begin{eqnarray*}
\frac{\partial \dot{V}}{\partial v}(s,t,x,v) & = &
-E_x(s,X(s,t,x,v))
\frac{\partial X}{\partial v}(s,t,x,v) \\
& = & - \rho(s,X(s,t,x,v)) \frac{\partial X}{\partial v}(s,t,x,v)
\end{eqnarray*}
so that
\begin{equation}
\label{vderiv} \frac{\partial V}{\partial v} (s,t,x,v) = 1 +
\int_s^t \rho(\tau, X(\tau,t,x,v)) \frac{\partial X}{\partial v}
(s,t,x,v) \ d\tau .
\end{equation}
Thus,
$$ \left \vert \frac{\partial X}{\partial v} (s,t,x,v) \right \vert \leq
\int_s^t \left \vert \frac{\partial V}{\partial v} (\tau,t,x,v)
\right \vert \ d\tau$$ and $$ \left \vert \frac{\partial
V}{\partial v} (s,t,x,v) \right \vert \leq 1 + \int_s^t \left
\vert \rho(\tau, X(\tau,t,x,v)) \frac{\partial X}{\partial v}
(s,t,x,v) \right \vert \ d\tau.$$ Combining the two inequalities,
we use Lemma $2$ and (\ref{rhobound}) so that
\begin{eqnarray*}
\left \vert \frac{\partial X }{\partial v} (s,t,x,v) \right \vert
+ \left \vert \frac{\partial V}{\partial v} (s,t,x,v) \right \vert
& \leq & 1 + \int_s^t\left ( \left \vert \frac{\partial
V}{\partial v} (\tau,t,x,v) \right \vert + C \left \vert
\frac{\partial X}{\partial v} (\tau,t,x,v) \right \vert \right ) \
d\tau \\
& \leq & 1 + C\int_s^t \left ( \left \vert \frac{\partial
V}{\partial v} (\tau,t,x,v) \right \vert + \left \vert
\frac{\partial X}{\partial v} (\tau,t,x,v) \right \vert \right ) \
d\tau
\end{eqnarray*}
for $\vert x \vert> R(t)$ and $\vert v \vert \leq Q_g(t)$. Then,
by Gronwall's Inequality, for $t \in [0,T]$, $s \in [0,t]$, $\vert
x \vert > R(t)$, and $\vert v \vert \leq Q_g(t)$
\begin{equation}
\label{charbound} \left \vert \frac{\partial X}{\partial v}
(s,t,x,v) \right \vert + \left \vert \frac{\partial V}{\partial v}
(s,t,x,v) \right \vert \leq C
\end{equation}
and we have bounds on $v$-derivatives of characteristics.\\

Now that $\rho(t,x)$ and $\frac{\partial X}{\partial v} (s,t,x,v)$
are bounded, define
$$ \Lambda(t)  := \sup_{\begin{array}{ccc} \vert x \vert > R(s) \\ s \in [0,t] \end{array}}
\vert \rho(s,x)\vert$$ and

$$ \Upsilon(t) :=
\sup_{\begin{array}{ccc} \vert x \vert > R(t) \\ \vert v \vert
\leq Q_g(t) \\ s \in [0,t] \end{array}} \left \vert \frac{\partial
X}{\partial v}(s,t,x,v) \right \vert.$$ Then, for $\vert v \vert
\leq Q_g(t)$ and $\vert x \vert > R(t)$,
\begin{eqnarray*}
\left \vert \int_0^t \rho(\tau, X(\tau,t,x,v)) \frac{\partial
X}{\partial v}(\tau,t,x,v) \ d\tau \right \vert & \leq & \int_0^t
\Lambda(\tau) \Upsilon(t) \ d\tau \\
& = & \Upsilon(t) \left ( \int_0^t \Lambda(\tau) \ d\tau \right ) \\
& =: & G(t).
\end{eqnarray*}
Notice $G \in C[0,T]$, $G$ increasing, $G(t) \geq 0$ for every $t
\in [0,T]$, and $G(0) = 0$. Define $$T_0 := \sup \left \{
\tilde{T} : G(\tilde{T}) \leq \frac{1}{2} \right \}. $$ Thus,
$G(T_0) \leq \frac{1}{2}$ and by the above computation,
$$\left \vert \int_0^t \rho(\tau, X(\tau,t,x,v)) \frac{\partial
X}{\partial v}(\tau,t,x,v) \ d\tau \right \vert \leq \frac{1}{2}$$
for every $t \leq T_0$. By (\ref{vderiv}), we find for $t \in
[0,T_0]$
\begin{eqnarray*}
\left \vert \frac{\partial V}{\partial v} (0) \right \vert & \geq
& 1 - \left \vert \int_0^t \rho(\tau, X(\tau,t,x,v))
\frac{\partial X}{\partial
v}(\tau,t,x,v) \ d\tau \right \vert \\
& \geq & 1 - \frac{1}{2} \\
& = & \frac{1}{2}
\end{eqnarray*}
so that
\begin{equation}
\label{lowerbound} \left \vert \frac{1}{\frac{\partial V}{\partial
v}(0)} \right \vert \leq 2
\end{equation} for $\vert v \vert \leq Q_g(t)$ and $\vert x \vert >
R(t)$. Once we show $\rho(t,x) = 0$ for $\vert x \vert > R(t)$
with $t \in [0,T_0]$, then $\Lambda(t) = G(t) = 0$ for every $t
\in [0,T_0]$, and it follows that $T_0 = T$, thus bounding $\vert
\frac{\partial
V}{\partial v}(0) \vert$ from below for all $t \in [0,T]$.\\

Let $\vert x \vert > R(t)$. Using (\ref{ftoF}), (\ref{vderiv}),
and (\ref{lowerbound}), we find

\begin{eqnarray*}
\rho(t,x) & = & \int(F(v) - F(V(0,t,x,v))) \ dv \\
& = & \int F(v) \ dv - \int F(w) \frac{1}{\frac{\partial
V}{\partial v}(0)} \ dw \\
& = & \int F(w) \left ( 1 - \frac{1}{\frac{\partial V}{\partial
v}(0)} \right ) \ dw \\
& = & \int F(w) \left ( \frac{\frac{\partial V}{\partial v}(0) -
1}
{\frac{\partial V}{\partial v}(0)} \right ) \ dw \\
& = & \int F(w) \left ( \int_0^t \rho(\tau, X(\tau,t,x,w))
\frac{\partial X}{\partial v}(\tau) \ d\tau \right ) \
\frac{1}{\frac{\partial V}{\partial v}(0)} \ dw.
\end{eqnarray*}
Thus, for $\vert x \vert > R(t)$, we use Assumption $(I)$,
(\ref{charbound}) and (\ref{lowerbound}) to find

\begin{eqnarray*}
\vert \rho(t,x) \vert & \leq & \int_{\vert w \vert \leq Q_g(t)}
\vert F(w)\vert \left ( \int_0^t \vert \rho(t,X(\tau)) \vert \
\left \vert \frac{\partial
X}{\partial v}(\tau) \right \vert \ d\tau \right ) \left \vert \frac{1}{\frac{\partial V}{\partial v}(0)} \right \vert \ dw \\
& \leq & C \int_{\vert w \vert \leq Q_g(t)} F(w) \ \int_0^t \vert
\rho(\tau, X(\tau)) \vert \ d\tau \ dw.
\end{eqnarray*}
Now, define $$\mathcal{P}(s) := \sup_{\vert x \vert > R(s)} \vert
\rho(s,x) \vert.$$ The above inequality becomes
$$ \mathcal{P}(t) \leq C \left ( \int F(w) \ dw \right ) \int_0^t
\mathcal{P}(\tau) \ d\tau.$$ By Gronwall's Inequality,
$$ \mathcal{P}(t) \leq 0.$$
Thus,
$$\mathcal{P}(t) = 0,$$
and
\begin{equation}
\label{rhocompact} \rho(t,x) = 0
\end{equation}
for $\vert x \vert > R(t)$.  As previously stated, since
(\ref{rhocompact}) holds for $t \in [0,T_0]$, we can conclude that
it does so for all $t \in [0,T]$.  Then, since $\rho(t,x)$ has
compact support, we use Lemma 2 and (\ref{rhocompact}) to find $$
\vert x \vert > R(t) \Rightarrow \rho(\tau, X(\tau,t,x,v)) = 0$$
for any $t \in [0,T]$, $\tau \in [0,t]$, $\vert v \vert \leq
Q_g(t)$. Using this with (\ref{vderiv}), we may conclude
$$ \frac{\partial V}{\partial v}(s,t,x,v) = 1 $$
for $\vert x \vert > R(t)$, $\vert v \vert \leq Q_g(t)$ and $s\in
[0,t]$. Thus, the proof of Theorem $1$ is complete.

\section*{Section 3}

In order to prove Theorem $2$, let us first define the current
density, $j : [0,T] \times \mathbb{R} \rightarrow \mathbb{R} $ by
$$ j(t,x) := \int v \ ( F(v) - f(t,x,v) ) \ dv.$$
A well known result, known as the equation of continuity, follows
from (\ref{VP}) :
\begin{equation}
\label{continuity}
\partial_t \rho(t,x) + \partial_x j(t,x) = 0.
\end{equation}
Using the second result of Theorem $1$, we find for $\vert x \vert
> R(t)$, $\vert v \vert \leq Q_g(t)$, and $s \in [0,t]$
$$V(s,t,x,v) = v + \gamma(s,t,x)$$ for some $\gamma$. In addition,
(\ref{char}) implies
$$V(s,t,x,v) = v + \int_s^t E(\tau, X(\tau, t, x, v)) \ d\tau.$$
Thus, for $\vert x \vert > R(t)$,
$$ E(\tau, X(\tau,t,x,v)) = E(\tau, X(\tau,t,x,0)). $$
Define for $\vert x \vert > R(t)$, $$\gamma(s,t,x) := \int_s^t
E(\tau, X(\tau, t, x, 0)) \ d\tau$$ so that
\begin{equation}
\label{vindep} V(0,t,x,v) = v + \gamma(0,t,x).
\end{equation}
Now, let $\vert x \vert > R(t)$.  By Assumption $(I)$,
(\ref{fchar}), and Lemma $2$, we find
\begin{eqnarray*}
j(t,x) & = & \int v \ (F(v) - f(t,x,v)) \ dv \\
& = & \int v \ F(v) \ dv - \int v \ f_0(X(0), V(0)) \ dv \\
& = & - \int v \ F(V(0))  \ dv \\
& = & - \int v \ F(v + \gamma(0,t,x)) \ dv \\
& = & - \int F(w) [w - \gamma(0,t,x)] \ dw \\
& = & \int F(w) \ \gamma(0,t,x) \ dw \\
& = & \left ( \int F(w) \ dw \right ) \int_0^t E(\tau,
X(\tau,t,x,0)) \ d\tau.
\end{eqnarray*}
Thus, from this relation, we find
\begin{equation}
\label{j} \frac{\partial}{\partial t} (j(t,x)) = \left ( \int F(w)
\ dw \right ) E(t,x).
\end{equation}\\
Since $\rho(t)$ has compact support, let us write supp$(\rho(t))
\subset [-L,L]$ for some $L > 0$, and so
$$ E(t,x) = \frac{1}{2} \left ( \int_{-L}^x \rho(t,y) \ dy - \int_x^L
\rho(t,y) \ dy \right ). $$ Also, notice
\begin{equation}
\label{Ebdy}
E(t,L) = -E(t,-L).
\end{equation}
Therefore, using (\ref{continuity})
\begin{eqnarray*}
\frac{\partial}{\partial t} (E(t,x)) & = & \frac{1}{2}
\left (\int_{-L}^x \rho_t(t,y) \ dy - \int_x^L \rho_t(t,y) \ dy \right) \\
& = & \frac{1}{2} \left ( - \int_{-L}^x j_x(t,y) \ dy + \int_x^L
j_x(t,y)
\ dy \right ) \\
& = & \frac{1}{2}[j(t,-L) - j(t,x) + j(t,L) - j(t,x)] \\
& = & \frac{1}{2}(j(t,-L) + j(t,L)) - j(t,x).
\end{eqnarray*}
Thus, using (\ref{j}) and (\ref{Ebdy}), we find for $\vert x \vert
> R(t)$
\begin{eqnarray}
\frac{\partial^2}{\partial t^2} (E(t,x)) & = &
\frac{1}{2}[j_t(t,-L) + j_t(t,L)] - j_t (t,x) \\
& = & \frac{1}{2} \left ( \int F(w) \ dw \right ) [E(t,-L) +
E(t,L)] - \left ( \int
F(w) \ dw \right ) E(t,x) \\
& = & - \left ( \int F(w) \ dw \right ) E(t,x) \label{Eode}
\end{eqnarray}
Now, for $\vert x \vert > R(t)$
\begin{eqnarray*}
E(t,x) & = & \frac{1}{2} sign(x) \int_{-L}^L \rho(t,y) \ dy \\
& = & sign(x) E(t,L).
\end{eqnarray*}
For all $t \in [0,T]$, define $$ e(t) := E(t,L),$$ and
$$ \omega := \left ( \int F(w) \ dw \right )^\frac{1}{2}.$$
Then, (\ref{Eode}) yields
$$ e^{\prime \prime}(t) = -\omega^2 e(t), $$
and so
\begin{equation}
\label{e}
e(t) = c_1 \sin(\omega t) + c_2 \cos(\omega t)
\end{equation}
for some $c_1, c_2 \in \mathbb{R}$. Thus, for $\vert x \vert >
R(t)$
\begin{equation}
\label{Eform}
E(t,x) = sign(x) (c_1 \sin(\omega t) + c_2
\cos(\omega t)).
\end{equation}

Now, we use (\ref{VP}) and Assumption $(I)$ to find for $\vert x
\vert
> R(t)$
\begin{eqnarray*}
E(0,x) & = & \frac{1}{2} sign(x) \int_{-R}^R \rho(0,y) \ dy \\
& = & \frac{1}{2} sign(x) \int_{-R}^R \int (F(v) - f_0(y,v)) \ dv \ dy \\
& =: & E^0 sign(x).
\end{eqnarray*}
However, by (\ref{Eform}) $$E(0,x) = (c_1 \sin(0) + c_2 \cos(0))
sign(x).$$ So, $c_2 = E^0$. Also, for $\vert x \vert > R(t)$,
$$ E_t(0,x) = \frac{1}{2} [j(0,-L) + j(0,L)] - j(0,x) = 0.$$
But, $E_t(0,x) = \omega c_1$, so that for non-trivial $E$, $c_1 =
0$. Finally, we may write
\begin{equation}
\label{Efinal} E(t,x) = E^0 sign(x) \cos (\omega t)
\end{equation} for
$\vert x \vert
> R(t)$ where
$$ E^0 = \frac{1}{2} \int_{-R}^R \int (F(v) - f_0(y,v)) \ dv \ dy $$
and
$$ \omega = \left (\int F(w) \ dw \right )^\frac{1}{2}.$$
So, the proof of the first part of Theorem $2$ is complete.

Now, we may use this result to show the second part of the
theorem.  First, by (\ref{char}), we know for $s \in [0,t]$, $$
V(s) = v + \int_s^t E(\tau, X(\tau)) \ d\tau.$$  Using
(\ref{Efinal}) and Lemma $2$ in this equation, for $\vert x \vert
> R(t)$,
\begin{eqnarray*}
V(s) & = & v + \int_s^t E^0 sign(x) \cos (\omega \tau) \ d\tau \\
& = & v + E^0 sign(x) \left ( \frac{1}{\omega} \right ) [ \sin
(\omega t) - \sin (\omega s) ].
\end{eqnarray*}
In addition, for $\vert x \vert > R(t)$
\begin{eqnarray*}
X(s) & = & x - \int_s^t V(\tau) \ d\tau \\
& = & x - sign(x) \int_s^t \left [ v + \frac{E^0}{\omega} (
\sin(\omega t)
- \sin(\omega \tau) ) \right ] \ d\tau \\
& = & x - sign(x) \left ( (t - s) \left [ v + \frac{E^0}{\omega}
\sin(\omega t)
\right ] + \frac{E^0}{\omega} \int_s^t \sin(\omega \tau) \ d\tau \right ) \\
& = & x - sign(x) \left ((t - s) \left [ v + \frac{E^0}{\omega}
\sin(\omega t) \right ] - \frac{E^0}{\omega} \left [ \cos(\omega
t) - \cos(\omega s) \right ] \right ) .
\end{eqnarray*}
Thus, for $\vert x \vert > R(t)$, we can explicitly calculate the
characteristics at $s = 0$ as
$$ X(0) = x - sign(x) \left ( t [ v + \frac{E^0}{\omega} \sin (\omega t) ] -
\frac{E^0}{\omega} [ \cos(\omega t) - 1 ] \right ) $$ and
$$ V(0) = v + \frac{E^0 sign(x)}{\omega} \sin(\omega t).$$
Therefore we use (\ref{fchar}) to find
$$ f(t,x,v) = f_0(X(0), V(0))$$ for $\vert x \vert > R(t)$ where
$X(0)$ and $V(0)$ are as above.  Finally, using $(\ref{IC})$ and
Lemma $2$, we conclude
$$ f(t,x,v) = F(V(0,t,x,v)) = F \left ( v + \frac{E^0 sign(x)}{\omega} \sin(\omega
t) \right )$$ for $\vert x \vert > R(t)$, and the proof of Theorem
$2$ is complete.  The final section will be devoted to the proofs of
Lemmas $1$ and $2$.

\section*{Section 4}
We complete the paper with the proofs of the lemmas.\\

\noindent \emph{\textbf{Proof of Lemma $1$}} : Let $T
> 0$ be given and $f$ be a solution of (\ref{VP}) on $[0,T]$.
Define for $t \in [0,T]$,
$$Q(t) : = \sup \{ \vert v \vert : \exists \ x \in \mathbb{R},
\tau \in [0,t] \ s.t. \ f(\tau,x,v) \neq 0 \}.$$

We show that all momenta characteristics are bounded as functions
of $s$. Using (\ref{star}) and (\ref{char}),
\begin{eqnarray*}
\vert V(0,t,x,v) \vert & = & \left \vert v + \int_0^t
E(s,X(s,t,x,v)) ds \right \vert \\
& \geq & \vert v \vert - C^{(1)}
\end{eqnarray*}
By definition of $Q(t)$, if $\vert V(0,t,x,v) \vert \geq Q(0)$, we
have for every $y \in \mathbb{R}$, $$f_0(y,V(0,t,x,v)) = 0.$$ But,
by the above equation, if $\vert v \vert \geq Q(0) + C^{(1)}$,
then
$$ \vert V(0,t,x,v) \vert \geq \vert v \vert - C^{(1)} \geq Q_g(0)$$
which implies that $f_0(y,V(0,t,x,v)) = 0$. So, if
$f_0(y,V(0,t,x,v)) \neq 0$, we must have
\begin{equation}
\label{vbound} \vert v \vert \leq Q(0) + C^{(1)}.
\end{equation}
Since we wish to consider only non-trivial $f$, (\ref{fchar})
implies that $f(X(0,t,x,v), V(0,t,x,v)) \neq 0$ for some $t \in
[0,T]$, $x,v \in \mathbb{R}$, and thus (\ref{vbound}) must hold.
Taking the supremum over $v$ of both sides in (\ref{vbound}), we
find
$$ Q(t) \leq Q(0) + C^{(1)} \leq C $$
for every $t \in [0,T]$.  Since $F$ is compactly supported, it
follows that $Q_g(t) \leq C$ for all $t \in [0,T]$, as well.  In
addition, for $t \in [0,T]$ and $\vert v \vert
> Q_g(t)$, we have
$$ f(t,x,v) = f_0(X(0),V(0)) = 0.$$

\noindent \emph{\textbf{Proof of  Lemma $2$}} : Let $T
> 0$ be given and $f$ be a solution of (\ref{VP}) on $[0,T]$.
Define $Q_g(t)$ for $t \in [0,T]$ as in (\ref{Qg}).  Consider
$\vert v \vert \leq Q_g(t)$ and use (\ref{char}) to obtain the
integral form of characteristics :
$$ V(s) = v + \int_s^t E(\tau,X(\tau,t,x,v)) \ d\tau $$
and
$$ X(s) = x - \int_s^t \left ( v + \int_\tau^t E(\bar{s},
X(\bar{s},t,x,v)) \ d\bar{s} \right ) \ d\tau.$$ Recall,
$$ R(t) := R + tQ_g(t) + \int_0^t \int_\tau^t \Vert E(\bar{s}) \Vert_\infty \
d\bar{s} d\tau. $$ Thus, for $\vert x \vert > R(t)$ and $s\in
[0,t]$,
\begin{eqnarray*}
\vert X(s,t,x,v) \vert & \geq & \vert x \vert - \vert v \vert
(t-s) - \int_s^t \int_\tau^t \vert E(\bar{s}, X(\bar{s},t,x,v))
\vert \ d\bar{s} \ d\tau \\
& \geq & \vert x \vert - Q_g(t) (t-s) - \int_s^t \int_\tau^t
\Vert E (\bar{s}) \Vert_\infty \ d\bar{s} \ d\tau \\
& > & R + t Q_g(t) - (t-s) Q_g(t) + \int_0^t \int_\tau^t \Vert
E(\bar{s}) \Vert_\infty \ d\bar{s} \ d\tau - \int_s^t \int_\tau^t
\Vert E(\bar{s}) \Vert_\infty \ d\bar{s} \ d\tau \\
& = & R + s Q_g(t) + \int_0^s \int_\tau^t \Vert E(\bar{s})
\Vert_\infty \ d\bar{s} \ d\tau \\
& \geq & R + s Q_g(s) + \int_0^s \int_\tau^s \Vert E(\bar{s})
\Vert_\infty \ d\bar{s} \ d\tau \\
& = & R(s)
\end{eqnarray*}
Thus, the proof of Lemma $2$ is complete.\\

\end{document}